\documentclass[12pt]{amsart}
\usepackage{amsmath,amsthm,amsfonts,amssymb}
\begin{document} 
\newcommand{\A}{{\mathbb A}}
\newcommand{\B}{{\mathbb B}}
\newcommand{\C}{{\mathbb C}}
\newcommand{\N}{{\mathbb N}}
\newcommand{\Q}{{\mathbb Q}}
\newcommand{\Z}{{\mathbb Z}}
\renewcommand{\P}{{\mathbb P}}
\renewcommand{\O}{{\mathcal O}}
\newcommand{\R}{{\mathbb R}}
\newcommand{\rc}{\subset}
\newcommand{\rank}{\mathop{rank}}
\newcommand{\trace}{\mathop{tr}}
\newcommand{\dimc}{\mathop{dim}_{\C}}
\newcommand{\Sing}{\mathop{Sing}}
\newcommand{\Spec}{\mathop{Spec}}
\newcommand{\Auto}{\mathop{{\rm Aut}_{\mathcal O}}}
\newcommand{\alg}[1]{{\mathbf #1}}
\newtheorem*{definition}{Definition}
\newtheorem*{claim}{Claim}
\newtheorem{corollary}{Corollary}
\newtheorem*{Conjecture}{Conjecture}
\newtheorem*{SpecAss}{Special Assumptions}
\newtheorem{example}{Example}
\newtheorem*{examples}{Examples}
\newtheorem*{remark}{Remark}
\newtheorem*{observation}{Observation}
\newtheorem*{fact}{Fact}
\newtheorem*{remarks}{Remarks}
\newtheorem{lemma}{Lemma}
\newtheorem{proposition}{Proposition}
\newtheorem{theorem}{Theorem}
\title{%
On a question of Koll\'ar.
}
\author {J\"org Winkelmann}
\begin{abstract}
We show: If a bounded domain in a Stein space covers a compact
complex space, it must be smooth. This give a negative answer
to a question of Koll\'ar.
\end{abstract}
%
\address{%
J\"org Winkelmann \\
Mathematisches Institut\\
Universit\"at Bayreuth\\
Universit\"atsstra\ss e 30\\
D-95447 Bayreuth\\
Germany\\
}
\email{jwinkel@member.ams.org\newline\indent{\itshape Webpage: }%
http://btm8x5.mat.uni-bayreuth.de/\~{ }winkelmann/
}
\thanks{
{\em Acknowledgement.}
The author was 
supported by 
the Mittag-Leffler-Institute and
the DFG Forschergruppe 790 ``Classification of algebraic
surfaces and compact complex manifolds''.
}

\maketitle
\section{Introduction}
The theory of the Kobayashi pseudodistance provides a link between
complex analysis and metric topology. We use this link to discuss
two topics.

The Shafarevich conjecture (see \cite{S})
 postulates that the universal covering
of a projective complex manifold ought to be holomorphically
convex. One important piece of evidence for this conjecture
is the the following result which was proved by Siegel
\cite{Si} in 1949:
{\em If the universal covering
of a complex compact manifold can be realized as a bounded domain in $\C^n$,
then this domain is holomorphically convex.} 

This result was later improved by Ivashkovich (\cite{I}) who showed:
{\em If the universal covering of a compact K\"ahler manifold can be realized
as an open subset $D$ in a complex manifold $M$, then $D$ is locally Stein,
i.e., for each point $p\in M$ there exists an open neighbourhood $W$ of
$p$ in $M$ such that $W\cap D$ is Stein or empty.}

Minimal model theory for projective manifolds suggests that one ought
to accept ``mild'' singularities. Hence it is natural to ask,
whether this result can be generalized to singular
spaces, and moreover to ask, as Koll\'ar did in \cite{K}, 
whether there exist compact spaces with a bounded singular domain
as universal covering.

We will prove that this  is impossible and we will also provide
an alternative proof for Siegels result.

Both statements will arise as corollaries of two more technical
results.

\section{The results}

We present the two more technical results, 
proposition~1 and 2, from which we deduce the main results.

\begin{definition}
The action of a group $G$ on a topological space $X$ is said
to be ``cocompact'', if there exists a compact subset $K\subset X$
such that $G\cdot K=X$.
\end{definition}

\begin{examples}
\begin{enumerate}
\item
If a group $G$ acts transitively on a  space $X$, then the action
is cocompact, because one point is compact.
\item
If $M$ is a compact (real or complex) manifold, then the action of
the fundamental group $\pi_1(M)$ on the universal covering $\tilde M$
by deck transformations is cocompact.
\item
Let $\Gamma$ be a discrete group acting properly discontinously on
a locally compact topological space $X$ with a compact quotient
$X/\Gamma$ and let $G$ be a group of homeomorphisms of $X$
containing $\Gamma$.
Then the $G$-action (as well as the $\Gamma$-action) on $X$
is cocompact.
\end{enumerate}
\end{examples}

\begin{proposition}
Let $G$ be a group acting cocompactly by isometries
on a locally compact metric topological space $X$.

Then $X$ is a complete metric space and for every
isometric embedding
$i:X\hookrightarrow Y$ of $X$
into a metric topological space $Y$ the image $i(X)$ is closed
in $Y$.
\end{proposition}
\begin{proof}
Let $K$ be a compact subset of $X$ with $G\cdot K=X$.

We define a function $\rho:X\to\R$ as follows: For every
$x\in X$ the value $\rho(x)$ is defined as the supremum of all
$r>0$ for which the closed ball
$\bar B_r(x)=\{p\in X:d(p,x)\le r\}$ is compact.
Note that $\rho(x)>0$ for every $x$ because $X$ is locally
compact.
Observe furthermore that 
\[
\bar B_{r-\epsilon}(y)\subset\bar B_r(x)
\]
if $d(x,y)\le\epsilon$.
This easily implies
\[
|\rho(x)-\rho(y)|\le d(x,y)\ \ \forall x,y\in X
\]
which in turn implies that $\rho$ is a continuous function.
Hence there is a minimum for $\rho$ on the compact set $K$.
Since $\rho$ is evidently invariant under all isometries of
$X$ and $G\cdot K=X$, it follows that there exists a constant
$c>0$ such that $\rho(x)>c$ for all $x\in X$.

Let $(x_n)$ be a Cauchy sequence in $X$. Then there exists a
natural number $N$ such that $d(x_n,x_m)< c$ if $n,m\ge N$.
Then $x_n\in \bar B_c(x_N)$ for all $n\ge N$.
But $\bar B_c(x_N)$ is compact, because $\rho(x_N)>c$.
Therefore every Cauchy sequence in $X$ is convergent, i.e.,
$X$ is complete.

To prove the second statement, assume the contrary. Then there
exists an isometric embedding $i:X\to Y$ and a sequence
$x_n$ in $i(X)$ which converges to a point in $Y\setminus i(X)$.
But this would imply that $(x_n)$ is a Cauchy sequence in $X$
which does not converge inside $X$ --- a contradiction to the
completeness of $X$.
\end{proof}

\begin{corollary}
A homogeneous locally compact metric space is complete.
\end{corollary}

\begin{corollary}
Let $X$ be a homogeneous complex manifold.

If $X$ is hyperbolic, the Kobayashi pseudodistance $d_X$ defines
a complete metric on $X$.

If the bounded holomorphic functions on $X$ separate the points,
the Caratheodory-pseudodistance $c_X$ defines a complete
metric on $X$.
\end{corollary}

\begin{corollary}
A homogeneous bounded domain in $\C^n$ is a complete metric
space with respect to both the Kobayashi and the
Caratheodory pseudometric.
\end{corollary}

\begin{theorem}
Let $D$ be a domain (=connected open submanifold)
in a hyperbolic Stein manifold $Z$.
Assume that there is a group $G$ acting on $D$ cocompactly by biholomorphic
transformations.

Then $D$ is Stein.
\end{theorem}
\begin{remark}
If $D$ is a bounded domain in $\C^n$, then $D$ is contained in some
ball $B_R=\{v\in\C^n:||v||<R\}$ which is a hyperbolic Stein manifold.

Similarily for a bounded domain
$D$ of an arbitrary Stein manifold $Z$: There is
an embedding $\zeta:Z\hookrightarrow\C^n$, and for $R>>0$
the preimage $\zeta^{-1}(B_R(0))=\{z\in Z:||\zeta(z)||<R\}$
is a hyperbolic Stein manifold containing $D$ as an open subset.
\end{remark}

\begin{proof}
Let $i:D\hookrightarrow E$ be the envelope of holomorphy
(see e.g.~\cite{GR}).
There is a natural projection $\pi:E\to Z$ which is locally biholomorphic.
It follows that $E$ is also hyperbolic (\cite{KH}, prop.~3.2.9).

Note that $E$ is a Stein space with $\O(E)=\O(D)$ (see e.eg.~\cite{GR}). 
Hence the points
of $E$ correspond to the closed maximal ideals in the ring of
holomorphic functions on $D$. Therefore each holomorphic automorphism 
of $D$ extends to an automorphism of $E$. 

Thus $G$ acts on $D$ by automorphisms which extend to $E$ and therefore
preserve the Kobayashi distance $d_E$ on $E$. Hence $D$ equipped with
the metric given by $d_E$ becomes a metric space on which $G$ acts
cocompactly and isometrically. It follows that $D$ is closed in $E$.
However, $D$ is always open in $E$, hence
$D$ being closed in $E$ implies that $D$ is a union of connected
components of $E$. Now $\O(E)\simeq\O(D)$ and $D$ being
connected imply that $E$ is connected.
Thus $D=E$. As a consequence, $D$ is Stein.
\end{proof}

\begin{corollary}[Siegel \cite{S}]
Let $M$ be a compact complex manifold with universal covering $D$.
If $D$ can be embedded into a Stein manifold $Z$ as a bounded domain, 
then $D$ is Stein.
\end{corollary}

\begin{corollary}
A homogeneous bounded domain in a Stein manifold is itself Stein.
\end{corollary}

(This is known due to the fundamental work of Gindikin, Pyatetski-Shapiro
and Vinberg on bounded homogeneous domains \cite{GPSV}).

\begin{proposition}\label{th2}
Let 
$D$ be a connected complex space, let $\Gamma$ be a group
of automorphisms of $D$ 
acting cocompactly on $D$ and let 
$Z$ denote a non-empty $\Gamma$-invariant subset of $D$.

Then every bounded holomorphic function on $D$ which vanishes
on $Z$ must be identically zero.
\end{proposition}

\begin{proof}
Fix $p\in Z$.
Let $f:D\to\C$ be a bounded holomorphic function vanishing on $Z$.
Assume that $f$ is not constant.
Without loss of generality we may assume that $\sup\{|f(z)|:z\in D\}=1$.
Let $\Delta=\{w\in\C:|w|<1\}$.
Then 
\[
\sup_{z\in D} d_\Delta(f(z),f(p))=+\infty
\]
since
\[
d_\Delta(w,0)=\log\frac{1+|w|}{1-|w|}\quad (w\in\Delta)
\]
(see \cite{KH}, p.21).

By the definition of the Caratheodory pseudodistance $c_D$ we have:
\[
c_D(q,z)\ge d_\Delta(0,f(z))
\]
for all $z\in D$, $q\in Z$.

It follows that
\[
\rho(z) \stackrel{def}{=} \inf_{q\in Z}c_D(q,z)
\]
is an unbounded continuous function on $D$.
But $\rho$ is $G$-invariant, because $Z$ is $G$-invariant and
\[
c_D(g(z),g(w))=c_D(z,w)\ \ \forall z,w\in D\forall g\in Aut(D).
\]
This leads to a contradiction: Every continuous function is bounded
on the compact set $K$ and therefore $G\cdot K=D$ implies that every
$G$-invariant continuous function on $D$ is bounded.

Hence there exists no such non-constant function $f$, i.e., if a bounded
holomorphic function $f$ on $D$ vanishes along $Z$, it must vanish 
identically on $D$.
\end{proof}

\begin{corollary}
Let $D$ be a bounded domain in an irreducible Stein space $X$, let $G$ be
a group acting cocompactly on $D$ and let $Z$ be a closed analytic
subset of $X$ containing a non-empty $G$-invariant subset of $D$.

Then $Z=X$.
\end{corollary}
\begin{proof}
Because $X$ is Stein, $Z$ can be defined by global holomorphic functions
on $X$. The restriction of a holomorphic function on $X$ to $D$ is
necessarily bounded, because $D$ is relatively compact in $X$.

Thus thm.~\ref{th2} implies that every holomorphic function vanishing
on $Z$ must also vanish on $D$. Hence $D\subset Z$, which implies $Z=X$,
because $X$ is irreducible.
\end{proof}

For ball quotients this specializes to the following fact:

\begin{corollary}
Let $X$ be a compact complex manifold with universal covering
 $\pi:B\to X$ 
 where $B=\{v\in\C^n:||v||<1\}$.

Let $Z$ be a non-empty closed analytic subset of $X$ and let $\epsilon>0$. 
Then $\pi^{-1}(Z)$
is not contained in any proper closed analytic subset of
\[
B_{1+\epsilon}=\{v\in\C^n : ||v||<1+\epsilon\}.
\]
\end{corollary}
\begin{proof}
This is immediate, because $B_{1+\epsilon}$ is a Stein manifold
containing $B$ as relatively compact connected open subset.
\end{proof}

\begin{theorem}
Let $X$ be a reduced Stein space, and let $D\subset$ be a bounded domain
on which a group $G$ acts cocompactly.

Then $D$ is smooth.
\end{theorem}

\begin{proof}
The singular locus $\Sing(X)$ of $X$ is a closed analytic subset of $X$.
Since $X$ is reduced, we have $\Sing(X)\ne X$.
Now every automorphism of $D$ must stabilize $\Sing(D)=D\cap\Sing(X)$.
Therefore the preceding corollary implies that $\Sing(D)=\{\}$, i.e.,
$D$ is smooth.
\end{proof}

In particular, we obtain a negative answer to the question of Koll\'ar
discussed in the introduction:
\begin{theorem}
Let $X$ be a reduced Stein space with a bounded domain $D$.
If there exists a group of automorphisms of $D$ acting properly
discontinuously on $D$ with compact quotient, then $D$ must be smooth.
\end{theorem}

\end{document}